\documentclass[12pt]{amsart}
\usepackage{amsmath, amssymb, amsfonts}
\usepackage{epsfig}

\newcommand{\pp}[1]{{\mathfrak #1}}

\newcommand{\length}{\ell}

\newcommand{\rank}{\operatorname{rank}}

\newcommand{\Hom}{\operatorname{Hom}}

\newcommand{\sym}{\mathcal{S}}

\newcommand{\br}{\operatorname{br}}
\newcommand{\e}{\operatorname{e}}
\newcommand{\fitt}{\operatorname{Fitt_0}}

\newcommand{\ree}{\mathcal{R}}

\newtheorem{lemma}{Lemma}[section]
\newtheorem{theorem}[lemma]{Theorem}
\newtheorem{corollary}[lemma]{Corollary}
\theoremstyle{definition}
\newtheorem{remark}[lemma]{Remark}
\newtheorem{assumption}[lemma]{Assumption}

\title{\sc Buchsbaum-Rim Multiplicities as Hilbert-Samuel Multiplicities\\}

\author{\sc C-Y. Jean Chan, Jung-Chen Liu, and Bernd Ulrich$^*$}
\date{ March~22, 2006\\
\indent Supported in part by the NSF  }

\begin{document}

\begin{abstract}
  We study the Buchsbaum-Rim multiplicity $\br(M)$ of a finitely
  generated module $M$ over a regular local ring $R$ of dimension~$2$
  with maximal ideal $\pp m$.  The module $M$ under consideration is
  of finite colength in a free module $F$.  Write $F/M \cong I/J$,
  where $J \subset I$ are $\pp m$-primary ideals of $R$.  We first
  investigate the colength $\ell(R/\pp a)$ of any $\pp m$-primary
  ideal $\pp a$ and its Hilbert-Samuel multiplicity $\e(\pp a)$ using
  linkage theory.  As applications, we establish several multiplicity
  formulas that express the Buchsbaum-Rim multiplicity of the module
  $M$ in terms of the Hilbert-Samuel multiplicities of ideals related
  to $I$, $J$ and an arbitrary general minimal reduction of $M$.  The
  motivation comes from E.~Jones' article~\cite{J} who applied
  graphical computations of the Hilbert-Samuel multiplicity to the
  Buchsbaum-Rim multiplicity.
\end{abstract}

\maketitle

Let $R$ be a local Cohen-Macaulay ring with maximal ideal $\pp m$ and
infinite residue field.  Let $\pp a$ be an $\pp m$-primary ideal. In
this paper, we study the connection between the {\em colength} of $\pp
a$, {\em i.e.}, the length $\ell( R/\pp a)$ of $R/\pp a$, and the
Hilbert-Samuel multiplicity $\e(\pp a)$ of $\pp a$. It is known for an
$\pp m$-primary ideal $\pp b$ contained in $\pp a$ that $\e(\pp a) =
\e(\pp b)$ if and only if $\pp b$ is a reduction of $\pp a$.
Furthermore, if $\pp b$ is a minimal reduction of $\pp a$, then
\begin{equation}\label{eq:ideal}  
\e(\pp a)=\e(\pp b)=\ell(R/\pp b). 
\end{equation}
However, $\e(\pp a)$ and $\ell(R/\pp a)$ are not equal in general. In
one of our main theorems, Theorem~\ref{rankone}, we express, under
certain conditions, the colength of $\pp a$ in terms of the
Hilbert-Samuel multiplicity of ideals which are in the same linkage
class of $\pp a$.

Equation (\ref{eq:ideal}) can be generalized to modules using the
Buchsbaum-Rim multiplicity of a module $M$, denoted $\br(M)$. Let $N
\subset M$ be submodules contained in a free module $F$ of finite rank
such that $\ell(F/N) < \infty$. It is known that $N$ and $M$ have the
same Buchsbaum-Rim multiplicity if and only if $N$ is a reduction of
$M$.  Similar to ideals, if $N$ is a minimal reduction of $M$, then
\begin{equation}\label{eq:module}
\br(M)=\br(N)= \ell(F/N)
\end{equation}
({\em cf.}  ~\cite{BR2},~\cite{KT},~\cite{SUV}).

In the case where $F$ has rank one, $M$ is an $\pp m$-primary ideal
and $\br(M)=\e(M)$. We see that the Buchsbaum-Rim multiplicity is a
generalization of the Hilbert-Samuel multiplicity to modules. Like the
Hilbert-Samuel multiplicity, it characterizes reductions.
Using the theory of reductions of modules, we
reduce the problem of finding formulas for the Buchsbaum-Rim
multiplicity to the relationship between
the colength and the Hilbert-Samuel multiplicity of ideals. 
The latter question is answered for arbitrary licci ideals in
Theorem~\ref{rankone}. As an application, we obtain formulas for the
Buchsbaum-Rim multiplicity of a two-dimensional module in terms of 
the Hilbert-Samuel multiplicities of a certain Fitting ideal and
ideals linked to it, see Theorem~\ref{shortformula}. We also prove
expressions for the Buchsbaum-Rim multiplicity that involve Bourbaki
ideals associated to the module, see
Theorems~\ref{shortformula2},~\ref{thmallrank}, and
Corollary~\ref{summary}. The last corollary contains the work of
\cite{J} as a special case.

The paper is arranged in the following way: Section~1 introduces the
notion of the Buchsbaum-Rim multiplicity and its basic properties. We
also include the definitions of some notation and theorems that will
be used in the later sections. In Section~2, we state and prove the
main theorem that relates the colength and the Hilbert-Samuel
multiplicity of $\pp m$-primary ideals in regular local rings of
dimension two. In Section 3, we discuss several multiplicity formulas
that express the Buchsbaum-Rim multiplicity of a module in terms of
the Hilbert-Samuel multiplicity of $\pp m$-primary ideals related to
the module. Jones~\cite{J} provides a method for computing the
Buchsbaum-Rim multiplicity of modules of a special type. In
Section~\ref{computation} we compare the multiplicity formulas
obtained in Section~\ref{generalization} to the results of \cite{J}.

\section{Introduction to the Buchsbaum-Rim Multiplicity}\label{intro}

In 1964, Buchsbaum and Rim~\cite{BR2} introduced and studied the
multiplicity that bears their names. It was further studied by
Gaffney, Kirby, Rees and many others, including Kleiman and Thorup who
investigated the geometric theory of the Buchsbaum-Rim multiplicity
in~\cite{KT}. In this paper, we study the connection between the
Buchsbaum-Rim multiplicity and the Hilbert-Samuel multiplicity.

Throughout the paper, we assume that $R$ is a Noetherian local ring of
dimension $d$ with maximal ideal $\pp m$. Let $\pp a$ be an $\pp
m$-primary ideal of $R$. There exists a polynomial $P_{\pp a}(n)$ of
degree $d$ such that $P_{\pp a}(n)= \length(R/{\pp a}^n)$ for large $n
\in \mathbb{N}$. This polynomial is called the Hilbert-Samuel
polynomial and the coefficient of $\frac{n^d}{d !}$ is the {\em
  Hilbert-Samuel multiplicity} $\e(\pp a)$.

The Buchsbaum-Rim multiplicity can be viewed as a generalization of
the Hilbert-Samuel multiplicity.  For a submodule $M$ of finite
colength in a free module $F$ of rank $r$, Buchsbaum and
Rim~\cite{BR2} prove that there exists a polynomial $\lambda(n)$ of
degree $d+r-1$ such that for all large $n \in \mathbb N$,
\[ \lambda(n) = \length( \sym_n(F)/ \ree _n(M)) ,\]
where $\sym(F) = \oplus _{n \geq 0} \sym_n(F)$ is the symmetric
algebra of $F$ and $\ree (M)= \oplus_{n \geq 0} \ree _n(M)$ is the
image of the natural map $\sym(M) \to \sym(F)$. Notice that the
algebra $\ree(M)$ is the $R$-subalgebra of $\sym(F)$ generated by $M$.
The {\em Buchsbaum-Rim multiplicity} $\br(M)$ is defined to be the
coefficient of $\frac{n^{d+r-1}}{(d+r-1)!}$ in the polynomial
$\lambda(n)$.  Buchsbaum and Rim showed that $\br(M)$ is a positive
integer if $M$ is a proper submodule of $F$. Notice that if $r=1$,
then $M$ is an $\pp m$-primary ideal of $R$, $\lambda(n) = P_M(n)$ and
$\br(M)=\e(M)$.

If $\operatorname{depth} R \geq 2$, then the inclusion $M \subset F$
where $\ell(F/M) < \infty$ can be identified with the natural
embedding of $M$ into its double dual $M^{**}$. Hence in this case
$\br(M)$ is independent of the embedding of $M$ into a free module.
Moreover, if $R$ is a two-dimensional regular local ring, one can
define the Buchsbaum-Rim multiplicity of {\em any} finitely generated
$R$-module $M$: simply consider the natural map from $M$ to $M^{**}$,
which is free in this case, and replace $M$ by its image under this
map.

Let $F$ be a free $R$-module of rank $r$, let $M$ be a submodule of
$F$ with $\length (F/M) < \infty$, and let $U$ be a submodule of $M$.
Again, we write $\ree(U)$ and $\ree(M)$ for the $R$-subalgebras of
$\sym(F)$ generated by $U$ and $M$, respectively.  We say that $U$ is
a {\it{reduction}} of $M$ if $\ree (M)$ is integral over $\ree (U)$ as
rings.  A minimal reduction of $M$ is a reduction that is minimal with
respect to inclusion.  Notice that a reduction $U$ of $M$ is minimal
if and only if its minimal number of generators is $r+d-1$, at
least when the residue field of $R$ is infinite~\cite{Rees1}.

On the other hand, if we fix a basis for $F$, then the submodule $M$
of $F$ is associated with a matrix, denoted by $\widetilde M$, whose
columns are the generators of $M$. Recall that the zeroth Fitting
ideal $\fitt(F/M)$ is the ideal generated by the maximal minors of
$\widetilde M$. This ideal is independent of the choices of the
generators of $M$ and the basis of $F$.
We recall a theorem by Rees relating reductions of ideals and modules:

\begin{theorem}$($\rm{Rees}\cite[1.2]{Rees1}$)$\label{Rees}
  The submodule $U$ of $M$ is a reduction of $M$ if and only if the subideal
  $\fitt(F/U)$ is a reduction of $\fitt(F/M)$.
\end{theorem}

If $U$ is a reduction of $M$, then $\br(U) = \br(M)$.  The converse
holds in the case $R$ is equidimensional and universal catenary ({\em
  cf.}~\cite{KT},~\cite{SUV}).

The following theorem relates the notions of reductions and 
Buchsbaum-Rim multiplicity.

\begin{theorem}$($\rm{Buchsbaum-Rim}\cite[4.5(2)]{BR2},
\rm{Kleiman-Thorup}\cite[5.3({\em i})]{KT},
\rm{Simis-Ulrich-Vasconcelos}\cite[2.5]{SUV},
\rm{Bruns-Vetter}\cite[2.8 and 2.10]{BV}$)$\label{thmBR}
  Assume that $R$ is a Cohen-Macaulay local ring with infinite residue
  field. If $U$ is a minimal reduction of $M$, then
\begin{equation*}
\br(M)=\br(U)=\length(F/U)=\length(R/\fitt(F/U)) .
\end{equation*}
\end{theorem}

\medskip

Let $R$ be a Cohen-Macaulay local ring of dimension $d \geq 2$ with
infinite residue field and let $M$ be a finitely generated module of
finite colength in a free module $F$ of rank $r$.
For such a module $M$, there exist ideals $J \subset I$ of height
$\geq 2$ such that $F/M$ is isomorphic to $I/J$.  In fact one can take
$I \cong F/G$ to be a Bourbaki ideal of $F$ with $G \subset M$ a free
submodule of rank $r-1$ and $J$ to be the image of $M$ in $I$ ({\em
  cf.}~\cite[Chapter~7 no.~4, Theorem~6]{B}, \cite[3.2(a)(c)]{SUV2}).  
Notice that if $d=2$ and $\mu(M) \leq 3$ then $I$ and
$J$ can be chosen to be complete intersections.  Since $M$ is its own
minimal reduction in this case, we obtain the following equalities by
Theorem~\ref{thmBR},
\begin{equation}\label{eqlength}
 \br(M)= \length(F/M) = \length(R/J) -\length(R/I)=\e(J)-\e(I).
\end{equation}
We see that the Buchsbaum-Rim multiplicity is connected to
the Hilbert-Samuel multiplicity in this special case ({\em cf.}~\cite{J}).
We are interested in such a relationship for arbitrary modules. By
Theorem~\ref{thmBR}, $\br(M)$ is equal to the colength of the Fitting
ideal corresponding to a minimal reduction of $M$. Thus, the question
can be reduced to investigating the connection between the colength and the
Hilbert-Samuel multiplicity of ideals.

\section{Colength and the Hilbert-Samuel multiplicity}\label{colength}

In a Cohen-Macaulay local ring $R$, two proper ideals $\pp a$ and $\pp
a_1$ are {\em linked} with respect to a complete intersection ideal
$\pp c$, denoted $\pp a \sim \pp a_1$, if
$\pp a_1=\pp c: \pp a$ and $\pp
a=\pp c: \pp a_1$. If $R$ is local Gorenstein and $\pp
a$ is unmixed of grade $g$ ({\em i.e.}, $\dim R_{\pp p}=g$ for all
associated prime ideals $\pp p$ of $R/\pp a$), it suffices to require $\pp
a_1=\pp c:\pp a$.  We say an ideal $\pp a$ is {\em in
  the linkage class of a complete intersection} (or $\pp a$ is {\em
  licci} for simplicity) if there are ideals $\pp a_1, \dots, \pp a_n$
with $\pp a \sim \pp a_1 \sim \cdots \sim \pp a_n$ and $\pp a_n$ a
complete intersection.

\begin{theorem}$($\rm{Huneke-Ulrich}~\cite[proof of 2.5]{HU1}$)$\label{huneke.ulrich.link}
  Let $(R, \pp m)$ be a Gorenstein local ring with infinite residue
  field and let $\pp a$ be a licci $\pp m$-primary ideal linked to a
  complete intersection in $n$ steps.  Then there exists a sequence of
  links $\pp a =\pp a_0\sim \pp a_1 \sim \cdots \sim \pp a_n$ such
  that $\pp a_n$ is a complete intersection, and $\pp a_i$ and $\pp
  a_{i+1}$ are linked with respect to a minimal reduction of $\pp
  a_i$.
\end{theorem}

\begin{theorem}\label{rankone}
  In the setting of Theorem~$\ref{huneke.ulrich.link}$, we have
\begin{equation*}
\length(R/\pp a) = \sum_{i=0}^{n} (-1)^i \e(\pp a_i).
\end{equation*}
\end{theorem}

\begin{proof} If $\pp a$ is a complete intersection, then
  $\length(R/\pp a)=\e(\pp a)$ and the assertion is clear. We assume
  that $\pp a$ is not a complete intersection and do induction on $n$.
  Let $\pp b_0$ be a minimal reduction of $\pp a$ such that $\pp
  a_1=\pp b_0:\pp a$. Notice that $\e(\pp b_0) = \e(\pp a)$. The
  quotient ring $R/\pp b_0$ is Gorenstein since $\pp b_0$ is generated
  by a regular sequence.  Moreover,
\[
(\pp b_0:\pp a)/\pp b_0 = \Hom_{R/\pp b_0}(R/\pp a,R/\pp b_0)=
\Hom_{R/\pp b_0}(R/\pp a,\omega_{R/\pp b_0})
  =D_{R/\pp b_0}(R/\pp a),
\]
where $\omega_{R/\pp b_0}$ is the canonical module of $R/\pp b_0$ and $D$
denotes the dualizing functor. Since the dualizing functor preserves
length, we have
\[ \length((\pp b_0:\pp a)/\pp b_0)=\length(R/\pp a). \]
Therefore
\[ \begin{array}{ccl}
\length(R/\pp a) &=& \length(R/\pp b_0)- \length(R/(\pp b_0:\pp a)) \\
             &=& \e(\pp b_0) - \length(R/\pp a_1)\\
             &=& \e(\pp a) - \length(R/\pp a_1). \end{array}\]
The result now follows by induction.
\end{proof}

\begin{corollary}\label{ingclosed} 
  Let $(R, \pp m)$ be a regular local ring of dimension $2$ with
  infinite residue field. If $\pp a$ is an integrally closed $\pp
  m$-primary ideal, then
\[ \ell( R/\pp a) = \sum_{i=1}^{\infty} (-1)^{i+1} \e({\rm Fitt}_i(\pp a)) . \]
\end{corollary} 

\begin{proof} 
  We induct on $\mu(\pp a)$, the case $\mu(\pp a)=2$ being obvious.
  Thus let $\mu(\pp a)= r+1 \geq 3$. Choose $\pp a_0, \dots, \pp a_r$ as in
  Theorem~\ref{rankone}, and notice that $\mu(\pp a_1) < \mu(\pp a)$.
  By Huneke and Swanson~\cite[3.4 and 3.5]{HS}, $\pp a_1$
  is integrally closed, $\pp a_1={\rm Fitt}_2(\pp a)$, and $\rm
  Fitt_i(\pp a_1) = \rm Fitt_{i+1}(\pp a)$ for every $i \geq 1$. Now
  we apply Theorem~\ref{rankone} and the induction hypothesis.
\end{proof} 

In Theorem~\ref{rankone} the colength of a licci ideal is expressed in
terms of Hilbert-Samuel multiplicities. This result applies to any
$\pp m$-primary perfect ideal in a two-dimensional Gorenstein local
ring with infinite residue field. For three-dimensional rings,
J.~Watanabe~\cite{JWat} has proved that every $\pp m$-primary perfect
Gorenstein ideal is licci.

\begin{theorem}\label{shortformula}
  Let $R$ be a Gorenstein local ring of dimension $2$ with infinite
  residue field, let $M$ be a submodule of finite colength in a free
  module $F$ of rank $r$, and let $U$ be a minimal reduction of $M$.
  \begin{itemize}
  \item[(a)] There exists a sequence of links $\fitt (F/U) = \pp a_0
  \sim \pp a_1 \sim \cdots \sim \pp a_{r-1}$ such that $\pp a_{r-1}$ is a
  complete intersection, and $\pp a_i$ and $\pp a_{i+1}$ are linked with
  respect to a minimal reduction of $\pp a_i$.
  \item[(b)]
\begin{equation*}
\br(M)=\e(\fitt(F/M)) + \sum_{i=1}^{r-1} (-1)^{i} \e(\pp a_{i}) .
\end{equation*}
  \end{itemize}
\end{theorem}

\begin{proof}
  To prove part (a), notice that $\pp a = \fitt(F/U)$ has height $2$
  and is generated by the maximal minors of an $r$ by $r+1$ matrix.
  Thus $\pp a$ can be linked to a complete intersection in $r-1$ steps
  $\pp a = \pp a_0 \sim \pp a_1 \sim \cdots \sim \pp a_{r-1}$
  (cf.~\cite{Ap},~\cite{AN},~\cite{G}).  By
  Theorem~\ref{huneke.ulrich.link} we may assume that $\pp a_i$ and
  $\pp a_{i+1}$ are linked with respect to a minimal reduction of $\pp
  a_i$.
  
  Part (b) follows from (a), Theorems~\ref{Rees}, \ref{thmBR} and
  \ref{rankone}.
\end{proof}

\begin{remark}\label{rmksec2}
  The ideals $\pp a_i$, $0 \leq i \leq r-1$, of
  Theorem~\ref{shortformula} can be obtained concretely in the
  following way: After applying general row and column operations to
  the matrix $\widetilde{M}$ presenting $F/M$, the ideal $\pp a_i$ is
  generated by the maximal minors of the matrix consisting of the
  last $r-2\lceil\frac{i-1}{2}\rceil$ rows and the last
  $r+1-2\lceil\frac{i}{2}\rceil$ columns of $\widetilde{M}$
  (\cite{AN}).
\end{remark}

\smallskip

As an immediate consequence of Corollary~\ref{ingclosed} and
Theorem~\ref{shortformula}, we obtain that if $\fitt(F/U)$ is
integrally closed, then
\[ \br(M) = \sum_{i=0}^{r-1} (-1)^i \e( {\rm Fitt}_i( F/U ) ). \]

The following remark provides another point of view on the formula in
Theorem~\ref{shortformula}. 
  
\begin{remark}\label{auslander} 
As in Remark~\ref{rmksec2} we apply general row and column operations
to the matrix $\widetilde M$, and then obtain an exact sequence 
\[ R^n \stackrel{\widetilde M}{\longrightarrow} F
\longrightarrow C_0=F/M \longrightarrow 0 . \]
The Auslander dual $D(C_0)$ of $C_0$ is presented by $\widetilde M
^*$, 
\[ F^* \stackrel{\widetilde M^*}{\longrightarrow} R^{n^*}
\longrightarrow D(C_0) \longrightarrow 0. \]
Let $C_1$ be the quotient of $D(C_0)$ modulo the submodule generated
by the image of the last $n-r+1$ basis elements of $R^{n^*}$. The
submatrix of $\widetilde M^*$ involving the top $r-1$ rows presents
$C_1$. 

Continuing this way, we obtain a
sequence of modules $C_0, \dots, C_{r-1}$, where $C_i$ is the quotient
of $D(C_{i-1})$ modulo the submodule generated by the last two
generators. 
Notice that $C_i$ is represented by the transpose of the matrix consisting
of the first $r-2\lceil\frac{i-1}{2}\rceil$ rows and the first
  $r+1-2\lceil\frac{i}{2}\rceil$ columns of $\widetilde{M}$ described
  in Remark~\ref{rmksec2}. Hence $\fitt(C_i) = \pp a_i$ and by
  Theorem~\ref{shortformula},
\[ \br(M) =  \sum _{i=0}^{r-1} (-1)^i \e(\fitt(C_i)) . \]

\end{remark}

\bigskip

\section{Multiplicity Formulas}\label{generalization}

In this section, we discuss other connections between the
Buchsbaum-Rim multiplicity of modules and the Hilbert-Samuel
multiplicity of ideals. In fact, we relate the Buchsbaum-Rim
multiplicity of $M$ to the Hilbert-Samuel multiplicity of a generic
Bourbaki ideal of $F$ with respect to $M$, see
Theorem~\ref{shortformula2}.  However, if there is a need to fix
a certain Bourbaki ideal $I$ of $F$, the result in
Theorem~\ref{shortformula2} does not apply anymore.  Instead
Theorem~\ref{thmallrank} takes care of these cases.

\begin{theorem}\label{shortformula2}
  Let $(R, \pp m)$ be a Gorenstein local ring of dimension $2$ with
  infinite residue field, let $M$ be a submodule of finite colength in
  a free module $F$ of rank $r$, let $U$ be a minimal reduction of
  $M$, and let $\pp a_i$ be ideals as in
  Theorem~$\ref{shortformula}(a)$.  Then there exists an $\pp
  m$-primary Bourbaki ideal $I$ of $F$ and a subideal $J \subset I$
  such that $F/M \cong I/J$ and
\[ \br(M) = \e(J) - \e(I) + \e(\pp a_2) +\cdots +(-1)^{r-1} \e(\pp a_{r-1}). \]
In particular, if $\rank M=2$, then there exist $\pp m$-primary ideals
$J \subset I$ such that $F/M \cong I/ J$ and
\[
\br(M) = \e(J) - \e(I).
\]
\end{theorem}

\begin{proof} 
  We may assume $r \geq 2$. Let $\pp b_0$ be a minimal reduction of
  $\pp a_0 =\fitt(F/U)$ defining the link $\pp a_0 \sim \pp a_1$. We
  can find generators $u_1, \dots, u_{r+1}$ of $U$ in $F$ so that 
 $\pp a_0$ and $\pp a_1$ are the ideals of maximal minors of the matrices
  $\widetilde U=(u_1 | \cdots | u_{r+1})$ and $\widetilde{V}= (u_1|
  \cdots | u_{r-1})$, and $\pp b_0$ is generated by the determinants
  of $(u_1| \cdots | u_{r-1}|u_{r+1})$ and $(u_1 | \cdots | u_{r-1}|u_r)$.

 Let $G$ be the submodule of $U$ generated by $u_1, \dots,
 u_{r-1}$. As $\pp a_1=I_{r-1}(\widetilde V)$ has height 2, it follows
 that $G$ is free and $\pp a_1 \cong F/G$ is an $\pp m$-primary
 Bourbaki ideal of $F$. Thus we may take $I$ to be $\pp a_1$. 
 
 Now let $J$ be the image of $M$ in $I$. Clearly $J\cong M/G$ and
 hence $I/J \cong F/M$. Notice that $\pp b_0$ is the image of $U$ in
 $I$. As $U$ is a reduction of $M$, it follows that $\pp b_0$ is a
 reduction of $J$. Since by definition $\pp b_0$ is also a reduction
 of $\pp a_0$, we deduce $\e(J)=\e(\pp b_0)=\e(\pp a_0)$. Now
 Theorem~\ref{shortformula} gives
\[ \br(M) = \e(J) - \e(I) + \e(\pp a_2) +\cdots +(-1)^{r-1} \e(\pp
a_{r-1}). \]

\vspace{-6mm}

\end{proof} 

We would like to point out that the result in
Theorem~\ref{shortformula2} does not hold for an arbitrary pair of
Bourbaki ideals $J \subset I$ satisfying $F/M \cong I/J$.  This case
is treated in our next result.  Theorem~\ref{thmallrank} provides an
expression for $\br(M)$ in terms of $\e(I)$ and $\e(J)$ if $I$ and $J$
are already specified. This is motivated by the work in Jones~\cite{J}
where it is necessary to choose $I$ and $J$ to be monomial ideals in
order to extend the graphical computation of the Hilbert-Samuel
multiplicity of monomial ideals to the Buchsbaum-Rim multiplicity of
modules. Jones also provides a class of examples where the formula of
Theorem~\ref{shortformula2} does not hold for arbitrary Bourbaki
ideals $J\subset I$.

\begin{assumption}\label{assume} Let $(R, \pp m)$ be a Gorenstein
local ring of dimension 2 with infinite residue field, let $M$ be a
submodule of finite colength in a free module $F$ of rank $r$, and
assume $M$ has no free direct summand. Write $F/M \cong I/J$, where
$J\subset I$ are $\pp m$-primary ideals, $I$ has finite projective
dimension, and $\mu(I) \leq r$. Since $M \subset \pp m F$, we have
$\mu(I/J) = \mu(F/M) = \mu(F) =r \geq \mu(I)$ and therefore $J\subset
\pp m I$. Thus the lift $F \longrightarrow I$ of the above isomorphism
is surjective by Nakayama's Lemma. It induces an isomorphism $I \cong
F/G$, where $G$ is a free submodule of $F$ of rank $r-1$. By
restriction we obtain $J \cong M/G$. 

  Let $s_1, \dots, s_{r-1}$ be generators of $G$ and let $z_r, \dots, z_{2r}$
  be generators of a minimal reduction $U$ of $M$. Thinking of $s_i
  \in F$ and $z_j \in F$ as column vectors we form the matrices
\[
\widetilde{L} = (s_1 | \cdots | s_{r-1}| z_r | \cdots |z_{2r}), \quad
 \widetilde{U} = ( z_r | \cdots |z_{2r}), \quad
 \widetilde{N} = ( s_1 | \cdots | s_{r-1}| z_{2r-1} | z_{2r} ) .  \]
 By performing row operations on $\widetilde{L}$ and by adding
 suitable linear combinations of columns of $\widetilde{L}$ to
 later columns we may achieve these properties: \begin{itemize} \item
 $s_1, \dots, s_{r-1}$ still generate $G$.

\item $z_r, \dots, z_{2r}$ still generate a minimal reduction $U$ of $M$.

\item the images of $z_{2r-1}, z_{2r}$ in $M/G = J$ generate a minimal
reduction $J'$ of $J$.

\item if for each $i$ with $0 \leq i \leq r-1$, $J_i$ denotes the ideal
of maximal minors of the matrix consisting of the last $r-2 \lceil
\frac{i-1}{2} \rceil$ rows and the last $r+1-2\lceil \frac{i}{2}
\rceil$ columns of $\widetilde{U}$, then $J_i$ and $J_{i+1}$ are linked
with respect to a minimal reduction of $J_i$ for $0 \leq i \leq r-2$.

\item if for each $i$ with $0 \leq i \leq r-1$, $J_i'$ denotes the
ideal of maximal minors of the matrix consisting of the last $r-2
\lceil \frac{i-1}{2} \rceil$ rows and the last $r+1-2\lceil
\frac{i}{2} \rceil$ columns of $\widetilde{N}$, then $J_i'$ and
$J_{i+1}'$ are linked with respect to a minimal reduction of $J_i'$
for $0 \leq i \leq r-2$.  
Notice that $J_{r-1} = J'_{r-1}$ and if $r$ is odd then also 
$J_{r-2}=J'_{r-2}$.
\end{itemize} 
Finally, let $I'$ be any minimal reduction of $I$, and
$\fitt(I/I')=I_0 \sim I_1 \sim \cdots \sim I_{r-3}$ be a sequence of
links as in Theorem~\ref{huneke.ulrich.link}. 
\end{assumption}

Note that for the last two conditions in \ref{assume}, one only has to
check that the two minors corresponding to the first two rows or
columns in the matrix of $J_i$ (or $J_i'$) generate a
reduction of $J_i$ (resp. $J_i'$). 
  
\begin{theorem}\label{thmallrank}
With assumptions as in $\ref{assume}$ one has
\begin{equation*}
\begin{array}[t]{ccccl}
\br(M) &=& \e(J) - \e(I) &+& (\e(\fitt(I/J)) +E_U) - (\e(\fitt(I/J'))+ E_N)  \\
       & & & + &(\e(\fitt(I/I')) + E_I),  \end{array}
\end{equation*}
where $E_U = \displaystyle{
 \sum_{i=1}^{2\lfloor \frac{r-2}{2}\rfloor } (-1)^i \e(J_i) }$,
$E_N = \displaystyle{
 \sum_{i=1}^{2\lfloor \frac{r-2}{2}\rfloor } (-1)^i \e(J_i') }$, and
$E_I = \displaystyle{ \sum_{i=1}^{r-3} (-1)^i \e(I_i) }$.
\end{theorem}

\begin{proof}
As $U$ is a reduction of $M$, Theorem~\ref{Rees} shows that
$I_r(\widetilde{U})$ is a reduction of $\fitt(F/M) = \fitt
(I/J)$. Therefore \[ \e(I_r(\widetilde U )) = \e( \fitt (I/J) ) .\]
Applying Theorem~\ref{rankone} to the ideals $I_r(\widetilde{U})$,
$I_r(\widetilde{N}) = \fitt(I/J')$ and $\fitt(I/I')$ we obtain
\begin{equation}\label{eq1}
 \mbox{\small $
 \length(R/I_r(\widetilde U)) = \left\{ \begin{array}{ll}
 \e( \fitt(I/J)) + E_U + (-1)^{r-1}\e(J_{r-1}) & \mbox{if $r$ is even} \\
 \e( \fitt(I/J)) + E_U + (-1)^{r-2}\e(J_{r-2}) + (-1)^{r-1}\e(J_{r-1}) &
  \mbox{if $r$ is odd } \end{array} \right. $ }   
\end{equation}
\begin{equation}\label{eq2}  
  \mbox{\small $ - \length(R/I_r(\widetilde N)) = \left\{ \begin{array}{ll} 
 - \e( \fitt(I/J')) - E_N - \e(J'_{r-1}) & \mbox{if $r$ is even} \\ 
  - \e( \fitt(I/J')) - E_N - (-1)^{r-2} \e(J'_{r-2}) - 
 (-1)^{r-1} \e(J'_{r-1})  & \mbox{if $r$ is odd} \end{array} \right.
  $ } \end{equation} 
\begin{equation}\label{eq3}
 \mbox{\small $ \length(R/\fitt(I/I')) = \e( \fitt(I/I')) + E_I. $ }
\end{equation}
Moreover by Theorem~\ref{thmBR},
\begin{equation*} \begin{array}{rcl}
\length(R/I_r(\widetilde{N})) - \length(R/\fitt(I/I'))
   & = & \length(I/J') - \length(I/I') \\  & = &
 (\length(R/J') - \length(R/I))- (\length(R/I')-\length(R/I)) \\
   & = & \e(J')-\e(I') \\ & = &  \e(J)-\e(I). 
 \end{array}
\end{equation*}
Thus, we have 
\begin{equation}\label{eqHSr}
\length(R/I_r(\widetilde{N})) - \length(R/\fitt(I/I')) = \e(J)-\e(I).
\end{equation}

Theorem~\ref{thmBR} also shows
\begin{equation}\label{eqMU}
\br(M)=\br(U)=\length(R/I_r(\widetilde{U})).
\end{equation}
Now, by adding equations (\ref{eq1}),
(\ref{eq2}), (\ref{eq3}), (\ref{eqHSr}) and applying (\ref{eqMU})
we obtain the multiplicity formula in Theorem~\ref{thmallrank}.
\end{proof}

We state the rank two and rank three cases as a corollary. The
multiplicity formulas have a more simple form in these cases. 

\begin{corollary}\label{summary}
  We use the assumption of $\ref{assume}$. 
\begin{itemize}
\item[(a)] If $r=2$ then 

\medskip

\noindent $ \br(M)=\e(J)- \e(I)  + \e(\fitt(I/J)) - \e(\fitt(I/J')) . $

\medskip
\item[(b)] If $r=3$ then
\begin{equation*}
\br(M)=\e(J) - \e(I)  + \e(\fitt(I/J)) - \e(\fitt(I/J')) + \e(\fitt(I/I')).
\end{equation*}

\end{itemize}
\end{corollary}

\begin{proof}
These results follow immediately from Theorem~\ref{thmallrank}.
If $r=2$, then the ideal $I$ is its own minimal reduction
and $\e(\fitt(I/I'))=0$.
\end{proof}

\medskip

\begin{remark}\label{counterex} 
 It should be pointed out that if a minimal reduction $J'$ of $J$ is
 general enough, then there exist minimal reductions $U$ of $M$
 such that Assumption~\ref{assume} is satisfied.
 The following example shows that the formula of
 Corollary~\ref{summary}(a) fails for a specific 
 $J'$, and therefore \ref{assume} does not hold for this $J'$. 

Let $R=k[x,y]_{(x,y)}$ and $M$ be a finitely generated module of 
finite colength in a free module $F$ of rank $2$ such that the presenting 
matrix of $F/M$ is 
\[ \left( \begin{array}{ccccc} 
 x^{16} & 0 & 0 & x^5y^4 & -y^{14} \\
 0 & y^{10} & x^8y^4 & 0 & x^{20} \end{array} \right). \]
Then $F/M \cong I/J$ where $I=(x^{20}, y^{14})$ and 
$J=(x^{36}, x^{25}y^4, x^8y^{18}, y^{24})$.  Note that 
$J'=(x^{36}+ y^{24}, x^{25}y^4)$ is a minimal reduction of $J$. 
The value on the right-hand side of the formula of 
Corollary~\ref{summary}(a) is 
\[ \e(J) -\e(I) + \e(\fitt(I/J)) - \e(\fitt(I/J')) = 
   744 - 280 + 546 - 594 = 416 \] 
while $\br(M) = 420$~(see \cite[page~50]{L} for details).  

This example also shows that $\e(\fitt(I/J')$ is not independent of
the choice of $J'$. 
\end{remark} 

\smallskip

\section{A Graphical Interpretation of the Buchsbaum-Rim Multiplicities}\label{computation}

In this section, we consider modules of rank two arising from monomial
ideals. We compare our formulas to the
result of E.~Jones~\cite{J}, who presented a graphical computation of
the Buchsbaum-Rim multiplicity in this case.
 
We assume $R=k[x,y]_{(x,y)}$ where $k$ is a field, and let $\pp m$
denote the maximal ideal of $R$.  Let $I$ and $J$ be $\pp m$-primary
monomial ideals with $J \subset \pp m I$, $\mu(I)=2$ and $\mu(J)\leq
3$.  Let $F$ be a free module of rank $2$ and $M$ a submodule of $F$
such that $M/F \cong I/J$. Jones computes the Buchsbaum-Rim
multiplicity of $M$ and shows that $\br(M) = \e(J) - \e(I)$ with a few
exceptions. For this one may assume that $k$ is infinite.

We write $I=(x^s, y^t)$ and may assume that $J=(x^{s+i}, x^dy^{t+e},
y^{t+j})$.  The module $M$ can be taken to be the image in $F=R^2$ of
the matrix 
\[ \widetilde{M} = \left( \begin{array}{cccc} 
 -y^t & x^i & 0 & 0 \\ x^s & 0 & x^d y^e & y^j \end{array} \right)
\]
In \cite{J} the modules $M$ are classified into seven
cases: In Figure~$c$, the point $T(s,t)$ corresponds to the monomial
$x^sy^t$ and similarly for other points including those in
Figures~$a$'s and $b$'s.

\smallskip
\begin{center}
 \begin{minipage}[t]{2in}
 \begin{center}
 \epsfig{file=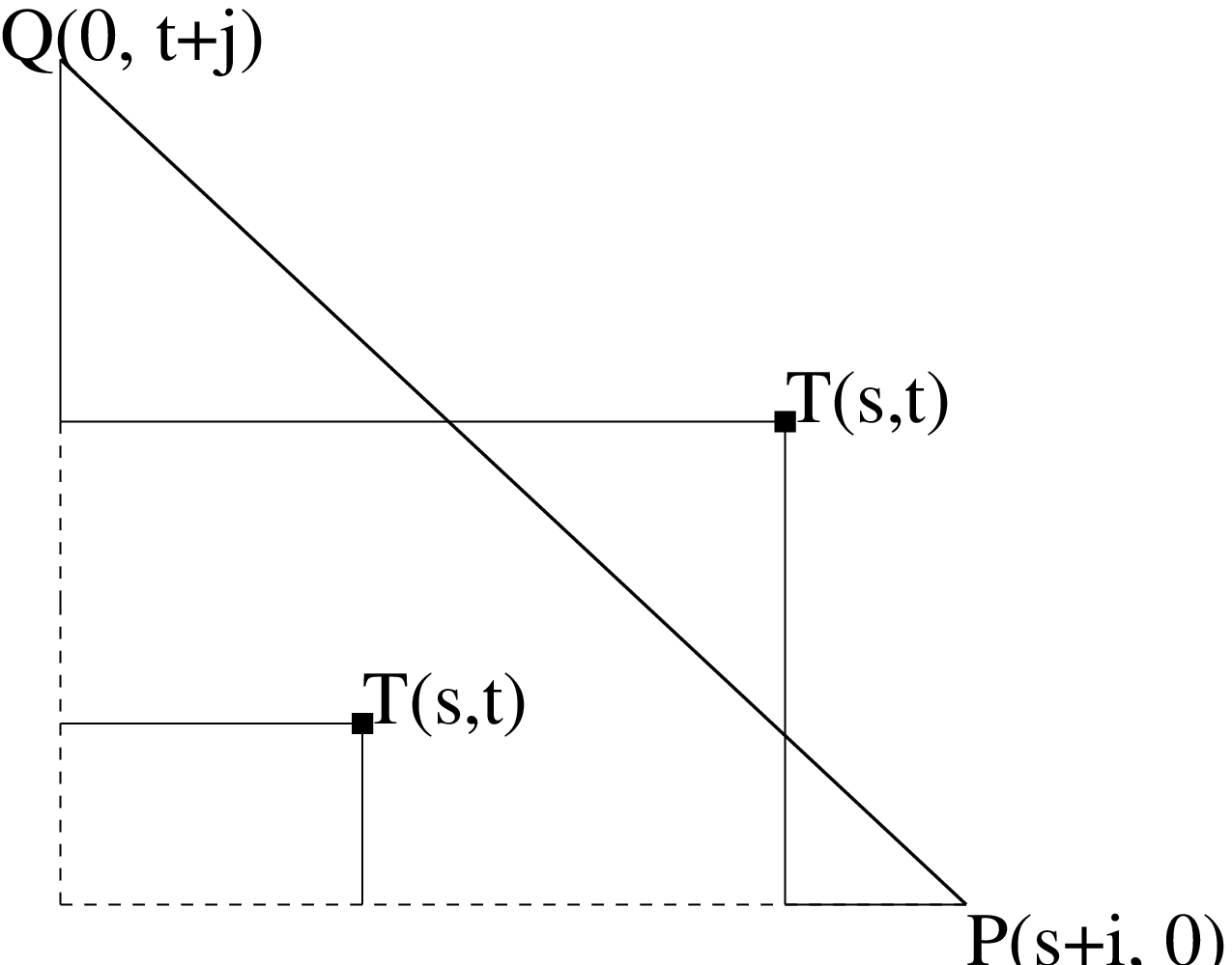, height=1.2in} \\
 {\small \em Figure~c. }
 \end{center}
 \end{minipage}
\end{center}
If $T$ is above the line segment $\overline{PQ}$ , then there are four
cases determined by the relative positions of the point $B(d, t+e)$
and $\overline{TQ}, \overline{PQ}, \overline{AQ}$ as shown in
Figures~$a$1--$a$4, where $\overline{AQ}$ is parallel to
$\overline{PT}$:

\smallskip

\begin{center}
 \begin{minipage}[t]{1in}
 \begin{center}
 \epsfig{file=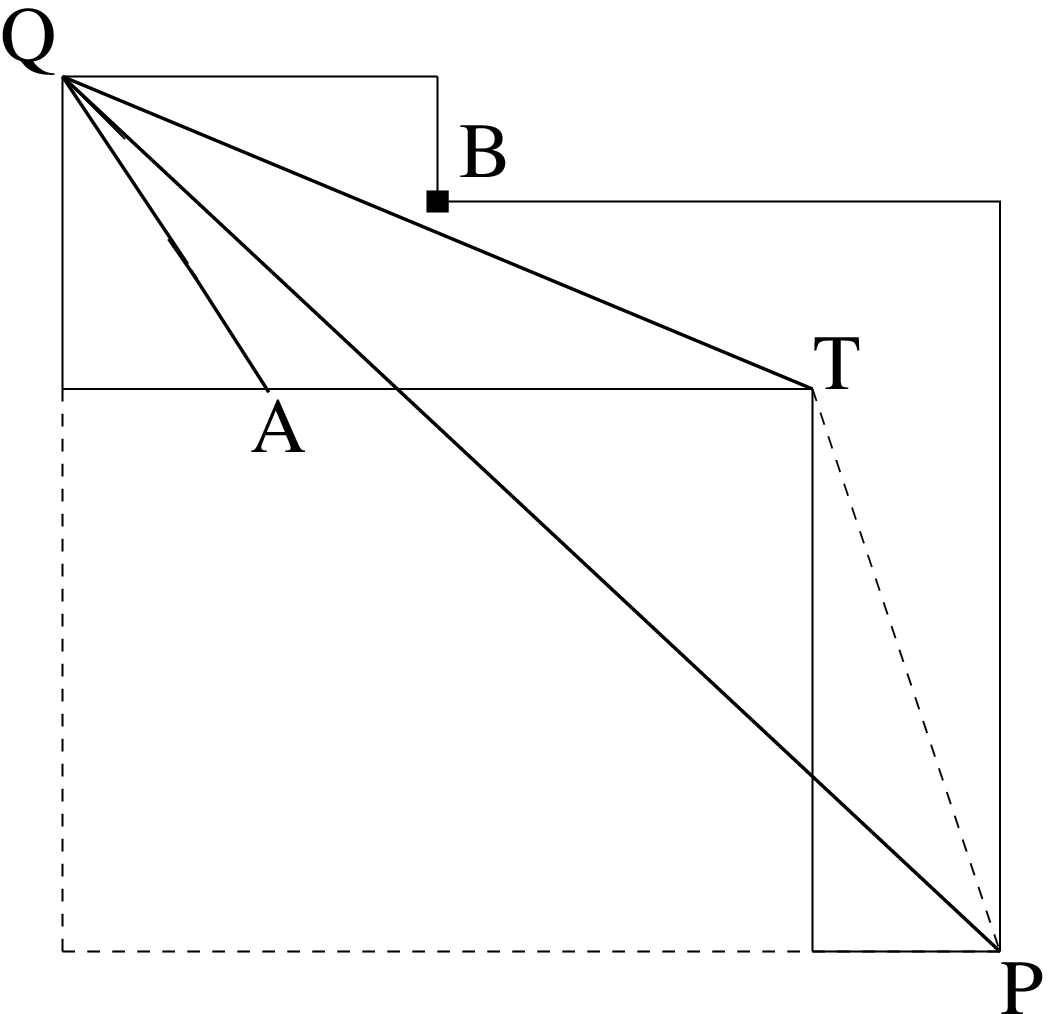, height=.9in} \\
 {\tiny \em Figure~a$1$. }
 \end{center}
 \end{minipage}
 \
 \
 \
 \begin{minipage}[t]{1in}
 \begin{center}
 \epsfig{file=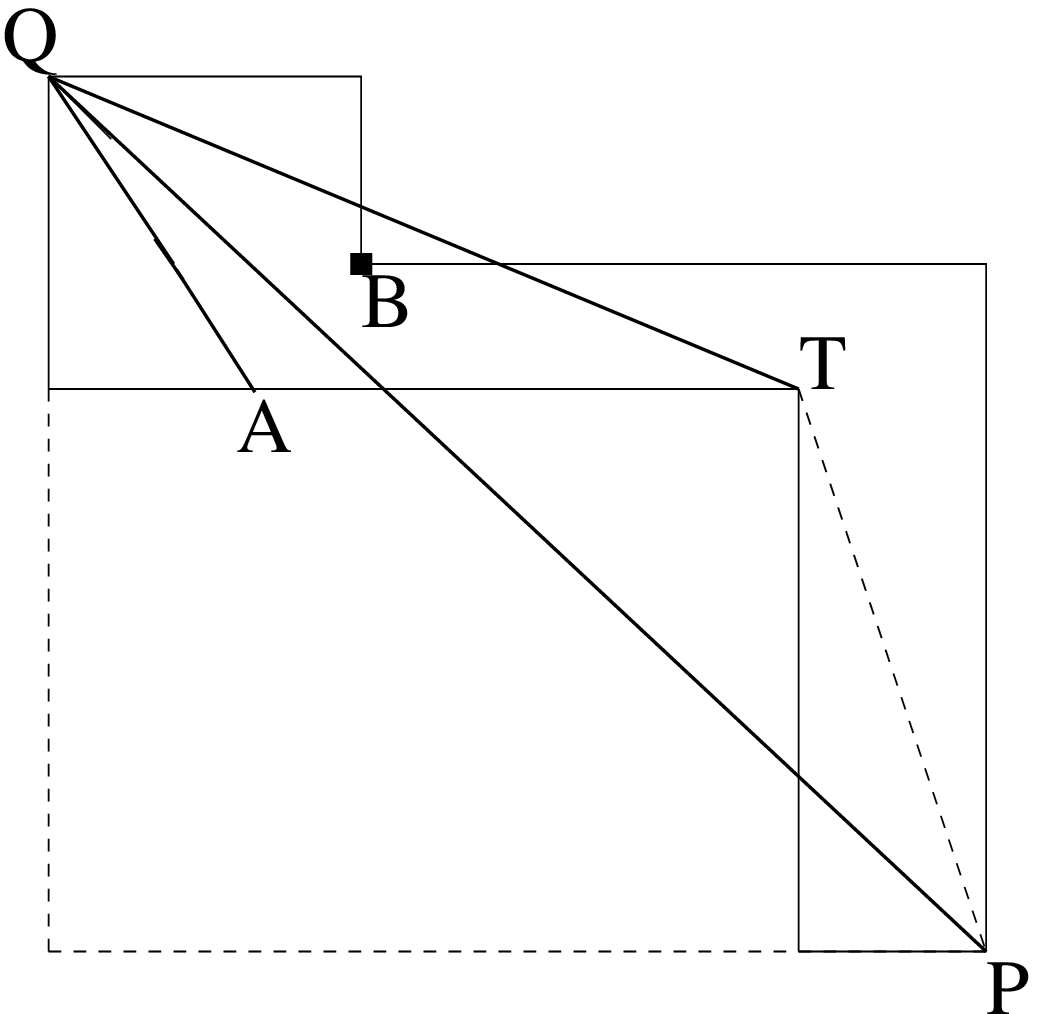, height=.91in} \\
 {\tiny \em Figure~a$2$. }
 \end{center}
 \end{minipage}
 \
 \
 \
\begin{minipage}[t]{1in}
 \begin{center}
 \epsfig{file=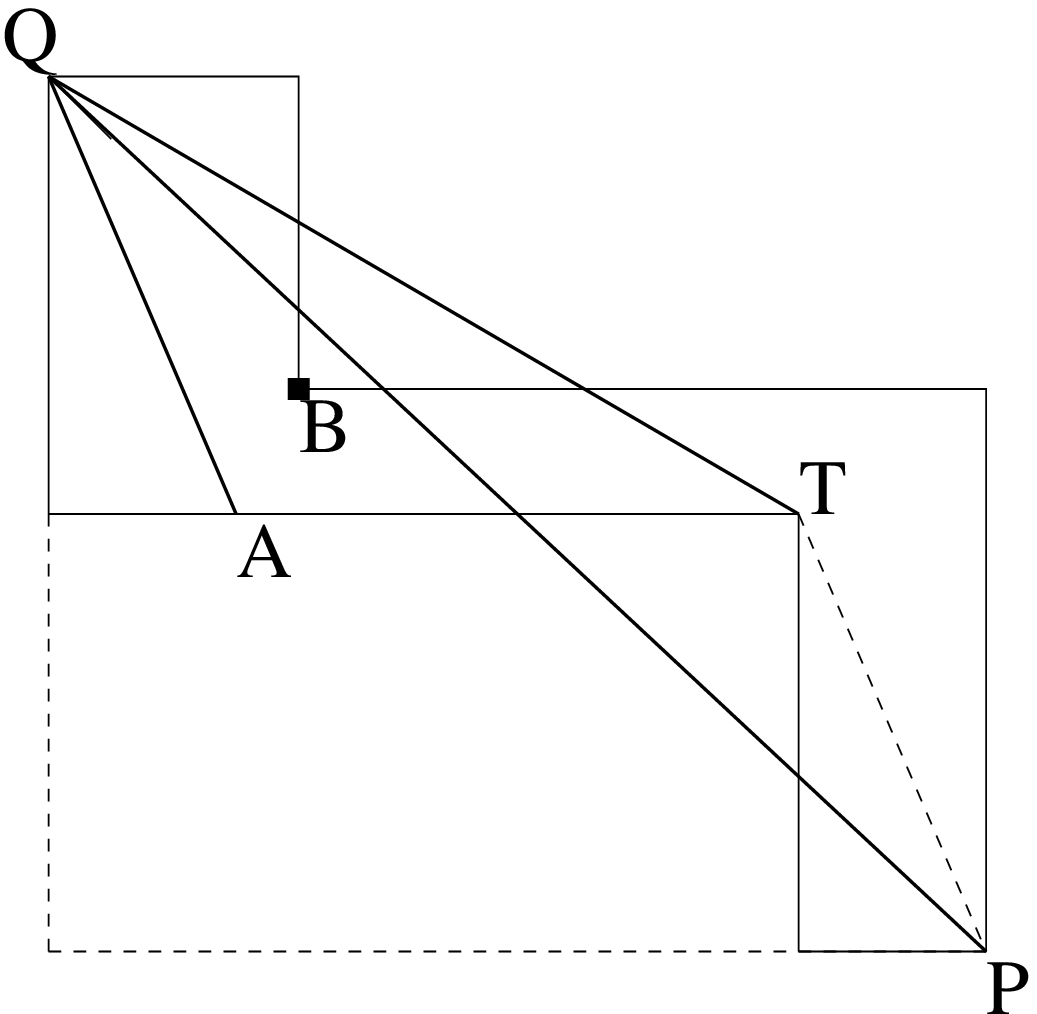, height=.91in} \\
 {\tiny \em Figure~a$3$. }
 \end{center}
 \end{minipage}
 \
 \
 \
 \begin{minipage}[t]{1in}
 \begin{center}
 \epsfig{file=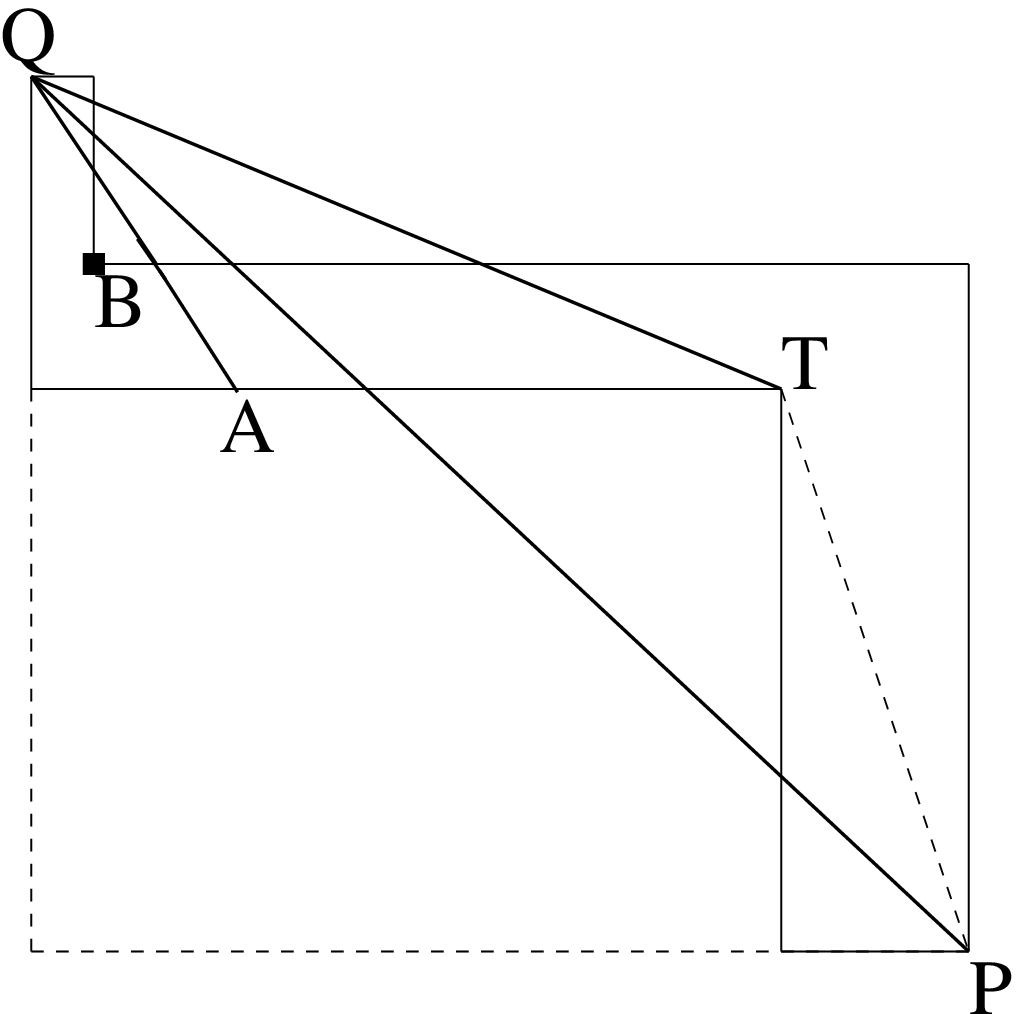, height=.91in} \\
 {\tiny \em Figure~a$4$. }
 \end{center}
 \end{minipage}
\end{center}
 If $T$ in Figure~$c$ is below
 $\overline{PQ}$, there are three cases determined by the relative
 positions of $B$ and $\overline{PQ}$, $\overline{PT}$ as shown in
 Figures~$b$1--$b$3:
\begin{center}
\begin{minipage}[t]{1in}
 \begin{center}
 \epsfig{file=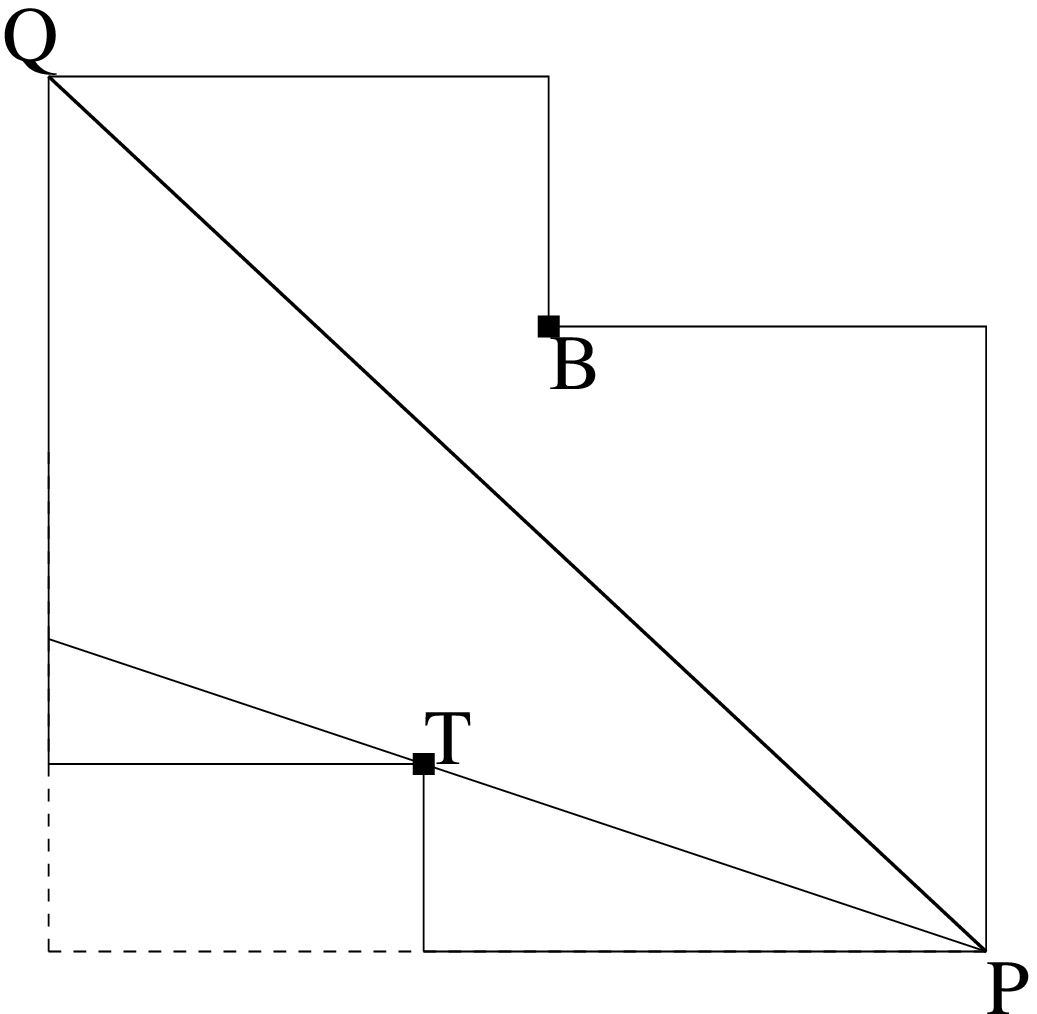, height=.9in} \\
 {\tiny \em Figure~b1. } \end{center} \end{minipage}
 \
 \
 \
 \
 \
 \begin{minipage}[t]{1in}
 \begin{center}
 \epsfig{file=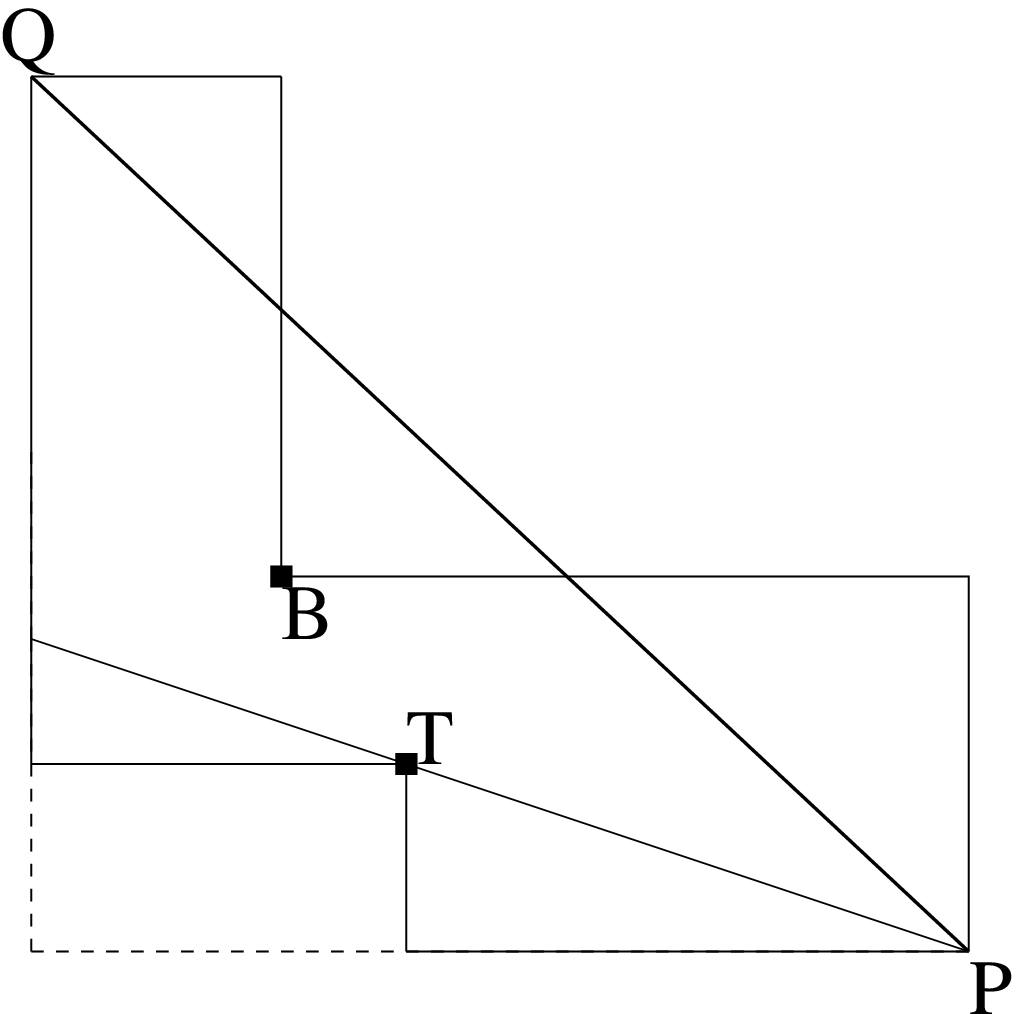, height=.9in} \\
 {\tiny \em Figure~b2. }
 \end{center}
 \end{minipage}
 \
 \
 \
 \
 \
 \begin{minipage}[t]{1in}
 \begin{center}
 \epsfig{file=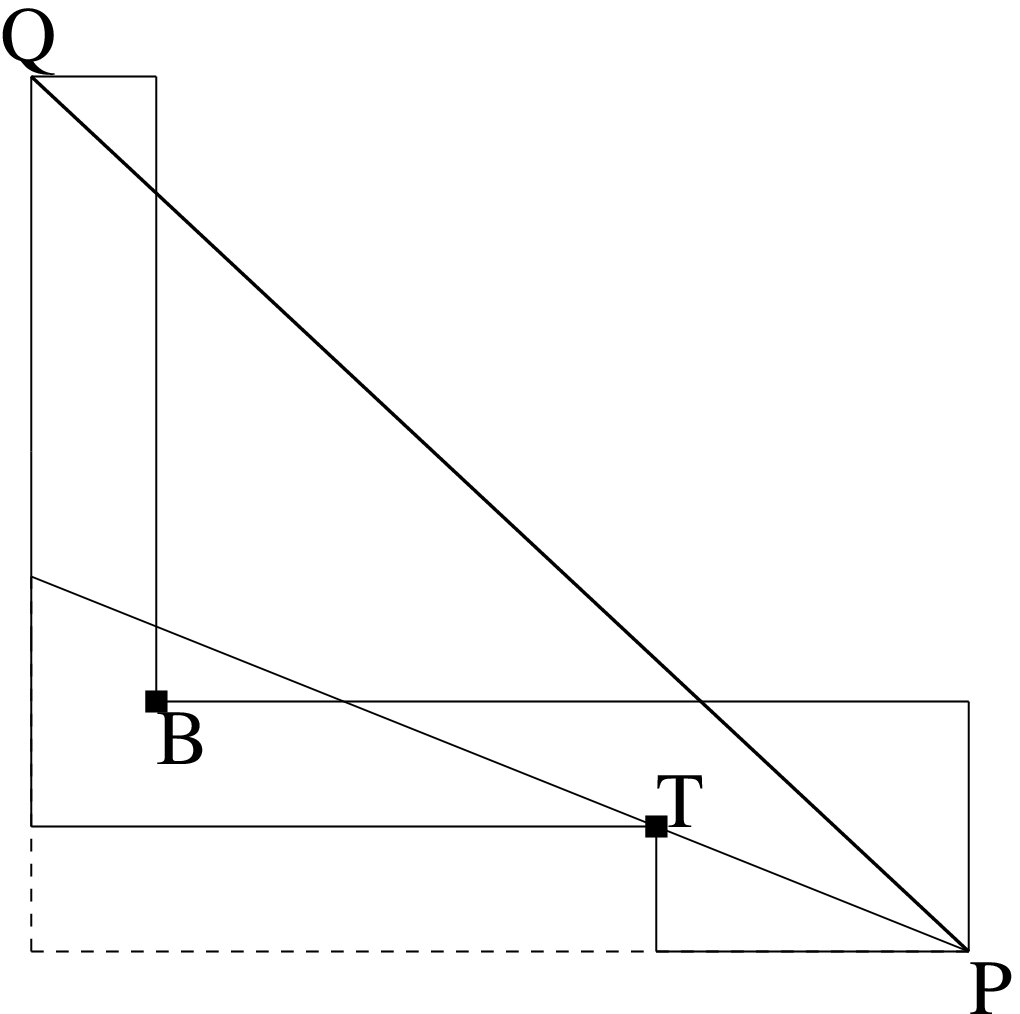, height=.9in} \\
 {\tiny \em Figure~b3. }
 \end{center}
 \end{minipage}
\end{center}

For the cases in Figures~$a1$ and $b1$, let $U$ be the submodule of
$F=R^2$ generated by the columns in the matrix
\[ \widetilde U =\left( 
\begin{array}{ccc}
  -y^t & x^i & 0 \\
  x^s & 0 & y^j \\
\end{array}
\right).
\]
Then $U$ is a minimal reduction of the module $M$.  Notice that the
first column in $\widetilde U$ is the syzygy of the ideal $I$ and the
image of $U$ in $J$ is a minimal reduction $J'$ of $J$.  Therefore in
\ref{assume}, we may take $\widetilde N$ to be $\widetilde U$ and
$\widetilde L$ to be $\widetilde U$ with the first column repeated. By
performing row operations on $\widetilde U$ and by adding suitable
linear combinations of columns of $\widetilde U$ to later columns we
have all the conditions required for Corollary~\ref{summary}. Since
$J'$ is the image of $U$ in $J$ and $U$ is a reduction of $M$,
$\fitt(I/J')$ is a reduction of $\fitt(I/J)$. Hence by
Corollary~\ref{summary}(a),
\[ \br(M)= \e(J) -\e(I). \]
This was also shown in in \cite{J}.
 
In Figures~$a4$, $b2$ and $b3$, let $U$ be the submodule of 
$F$ generated by the columns in the matrix
\[ \widetilde{U} = \left(
\begin{array}{ccc}
  -y^t & x^i & 0 \\
  x^s & y^j & x^dy^e \\
\end{array}
\right).
\]
By the same argument, $\br(M)= \e(J) -\e(I)$. 

For the remaining cases, the modules of Figures~$a2$ and $a3$, we use the
computation of the Buchsbaum-Rim multiplicity given in \cite{J}. There it is 
shown that $M$ is a reduction of the module generated by $M$ and the 
vector $(0, x^s)$ in $F$, which is a direct sum of two monomial ideals. 
This allows for a computation of $\br(M)$. Thus in the case of Figure~$a2$, 
\begin{equation}\label{graph1}
\br(M)= \e(J) - \e(I) - 2 \cdot \mbox{dark area},
\end{equation}
where the dark area is the area of the triangle $TBQ$
indicated in the following Figure~$a2'$. On the other hand, the modules 
of Figure~$a3$ have
Buchsbaum-Rim multiplicity 
\begin{equation}\label{graph2} \br(M)=\e(J)-\e(I) -2\cdot\mbox{dark
area} + 2\cdot\mbox{light area},  \end{equation}
where the dark area is the area of the triangle $TBQ$
and the light area is the area of the triangle $PBQ$ as indicated in
Figure~$a3'$. 

\smallskip

\begin{center}
\begin{minipage}[t]{2in}
\begin{center}
\epsfig{file=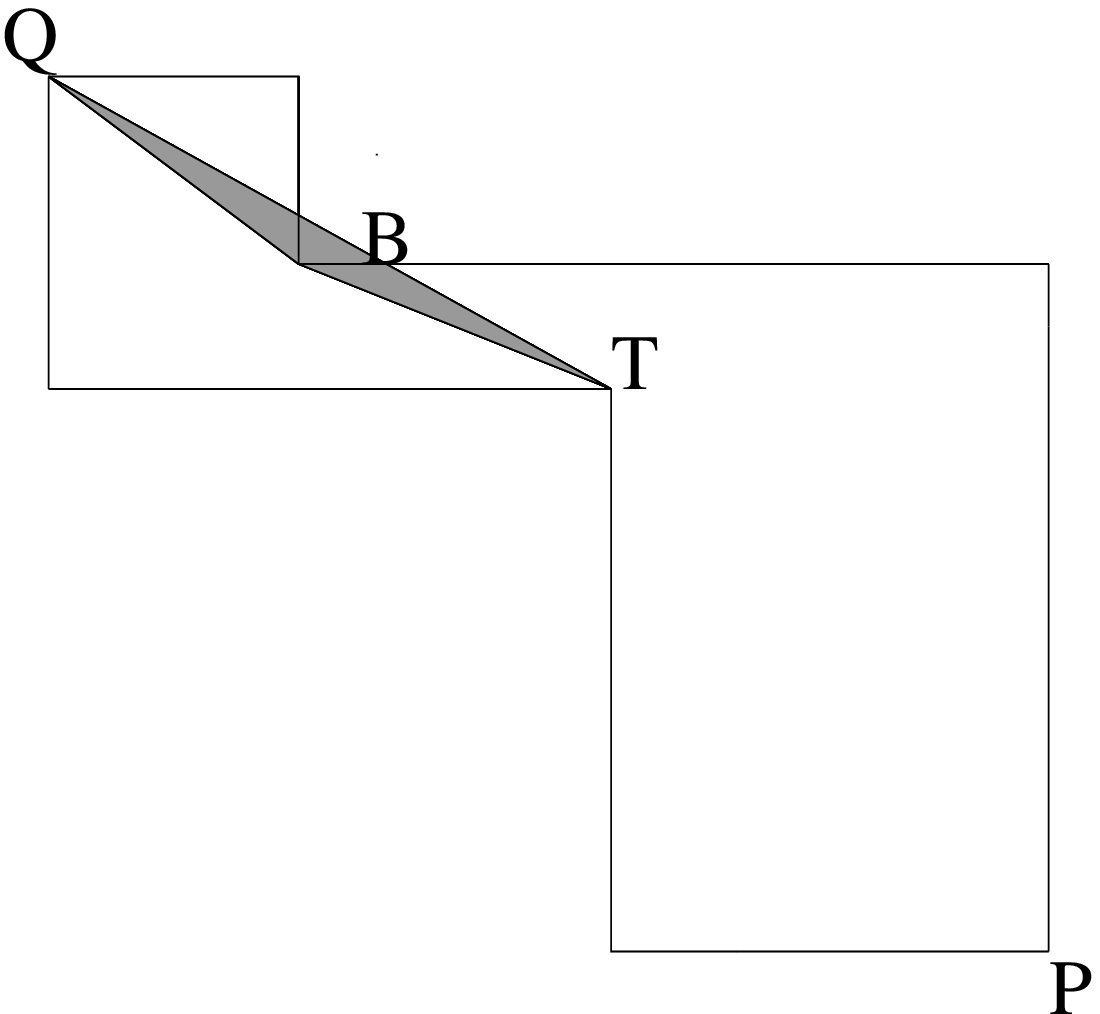, height=1.5in} \\
{\small \em Figure~a$2'$. } \end{center}
\end{minipage}
\
\begin{minipage}[t]{2in}
\begin{center}
\epsfig{file=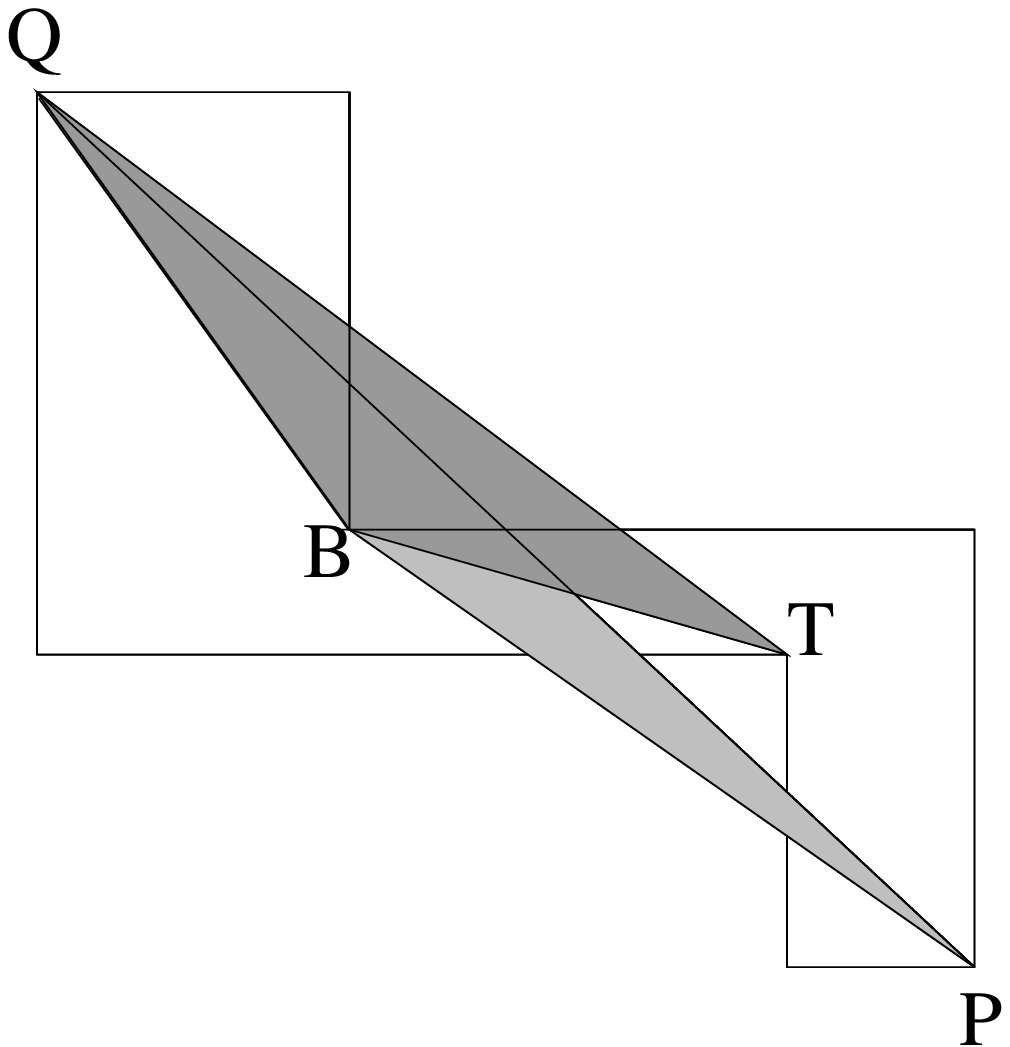,height=1.5in} \\
{\small \em Figure~a$3'$.} \end{center}
\end{minipage}
\end{center}

\smallskip

\noindent By Corollary~\ref{summary}(a), the extra terms subtracted 
in (\ref{graph1}) and (\ref{graph2}) are exactly 
\[ \e(\fitt(I/J')) - \e( \fitt(I/J) ) \]
for some general minimal reduction $J'$ of $J$ that satisfies
\ref{assume}.  We remark that in the first five cases, since
$\fitt(I/J)$ has a simple form, one can find a minimal
reduction $U$ of $M$ that is close to being monomial. For the cases
$a2$ and $a3$, this is much more complicated. 

\bigskip

\begin{center} 
{\sc Acknowledgment. } \\
The authors would like to thank Liz Jones for many valuable 
discussions. 
\end{center}

\bigskip
{\small
\noindent Department of Mathematics, University of Arkansas, Fayetteville,
AR~72701, USA \\
e-mail: cchan@uark.edu

\smallskip
\noindent 
Department of Mathematics, National Taiwan Normal University, Taipei,
Taiwan \\
e-mail: liujc@math.ntnu.edu.tw 

\smallskip
\noindent 
Department of Mathematics, Purdue University, West Lafayette,
IN~47907, USA \\
e-mails: ulrich@math.purdue.edu }

\end{document}